\documentclass{mcom-l}

\usepackage[cmtip,all]{xy}
\usepackage{amssymb}
\usepackage{algorithm}
\usepackage{algorithmic}

\newtheorem{theorem}{Theorem}[section]
\newtheorem{lemma}[theorem]{Lemma}

\theoremstyle{definition}
\newtheorem{definition}[theorem]{Definition}

\newtheorem{condition}[theorem]{Genericity condition}

\theoremstyle{remark}
\newtheorem{remark}[theorem]{Remark}

\numberwithin{equation}{section}

\DeclareMathOperator{\lcm}{lcm}

\newcommand{\al}{\alpha}
\newcommand{\be}{\beta}
\newcommand{\ga}{\gamma}
\newcommand{\de}{\delta}
\newcommand{\ep}{\epsilon}
\newcommand{\la}{\lambda}

\newcommand{\Z}{\mathbb{Z}}
\newcommand{\C}{\mathbb{C}}
\newcommand{\Q}{\mathbb{Q}}
\newcommand{\Fp}{\mathbb{F}_{p}}
\newcommand{\PP}{\mathbb{P}}

\newcommand{\abcd}{\left(\al:\be:\ga:\de\right)}
\newcommand{\XYZT}{\left(X:Y:Z:T\right)}
\newcommand{\K}{\mathcal{K}}
\newcommand{\mO}{\mathcal{O}}
\newcommand{\mC}{\mathcal{C}}
\newcommand{\mCt}{\tilde{\mathcal{C}}}
\newcommand{\mE}{\mathcal{E}}
\newcommand{\J}{Jac\left(\mathcal{C}\right)} 
\newcommand{\Jt}{Jac\left(\tilde{\mathcal{C}}\right)} 

\newcommand{\tun}{\theta_{1}^2 }
\newcommand{\tde}{\theta_{2}^2 }
\newcommand{\ttr}{\theta_{3}^2 }
\newcommand{\tqu}{\theta_{4}^2 }
\newcommand{\tci}{\theta_{5}^2 }
\newcommand{\tsi}{\theta_{6}^2 }
\newcommand{\tse}{\theta_{7}^2 }
\newcommand{\thu}{\theta_{8}^2 }
\newcommand{\tne}{\theta_{9}^2 }
\newcommand{\tdi}{\theta_{10}^2 }
\newcommand{\ton}{\theta_{11}^2 }
\newcommand{\tdo}{\theta_{12}^2 }
\newcommand{\ttz}{\theta_{13}^2 }
\newcommand{\tqo}{\theta_{14}^2 }

\newcommand{\tsz}{\theta_{16}^2 }
\newcommand{\zero}{\left(0\right)}
\newcommand{\z}{\left(z\right)}

\begin{document}

\title{Factorization with genus 2 curves}


\author{Romain Cosset}
\address{LORIA\\ Campus Scientifique - BP 239 \\ 54506 Vandoeuvre-l\`es-Nancy , France}
\curraddr{}
\email{romain.cosset@loria.fr}

\subjclass[2000]{Primary 11Y05; Secondary 11Y16, 11Y40}

\date{}

\dedicatory{}

\begin{abstract}
The elliptic curve method (ECM) is one of the best factorization methods available. It is possible to use hyperelliptic curves instead of elliptic curves but it is in theory slower. We use special hyperelliptic curves and Kummer surfaces to reduce the complexity of the algorithm. Our implementation GMP-HECM is faster than GMP-ECM for factoring big numbers.
\end{abstract}

\maketitle

\section*{Introduction}

The elliptic curve method (ECM) introduced in 1985 by H. W. Lenstra, Jr. \cite{Lenstra87} plays an important role in factoring integers. ECM is used to find ``medium sized'' (up to $60$ digits) prime factors of ``random'' numbers. It is also used by other factoring algorithms like sieving methods for cofactoring \cite{Kruppa08}. ECM is a generalization of Pollard's $p-1$ algorithm: instead of working in $\Fp^{*}$, we work in the group of points on an elliptic curve. We generalize it by using hyperelliptic curves of genus~$2$ instead of elliptic curves of genus $1$.

\bigskip

At first sight the ``hyperelliptic curve method'' (HECM) seems slower that ECM because of two reasons: first the arithmetic of hyperelliptic curves is slower compared to the arithmetic of elliptic curves, and secondly its probability of success is smaller than that of ECM.\\
Indeed this probability depends on the probability of the cardinality of the Jacobian of the curve modulo a prime factor $p$ of $n$ being smooth. But the probability of a number being smooth decreases with its size. By Weil's theorem, the cardinality of a Jacobian of a genus $2$ curve on $\Fp$ is around $p^2$ whereas the cardinality of an elliptic curve over the same field is around $p$. This seems like a major deterrent to using HECM. To avoid it we use special hyperelliptic curves whose Jacobians are isogenous to the product of two elliptic curves: they are called decomposable. One run of HECM with this kind of curves is equivalent to two simultaneous runs of ECM, thus the probability of success of HECM with decomposable hyperelliptic curves is comparable to that of two ECM. \\
The complexity of stage $1$ is dominated by the cost of the arithmetic on the curves. In genus $1$, the quickest arithmetic is obtained by using Montgomery formul\ae{} \cite{Montgomery87} or the new Edwards curves \cite{BeBiLaPe08}. In genus $2$, the best known formul\ae{} use Kummer surfaces \cite{Gaudry07}. The Kummer surface is a variety obtained by identifying opposite points of the Jacobian. However, not all genus $2$ hyperelliptic curves map to a rational Kummer surface.

\bigskip

The curves used in HECM are built as reduction modulo $n$ of curves defined over $\Q$ which are selected for satisfying different constraints. We need a large number of possible trial curves in order to factor $n$, hence the requirement that there exists an infinite family of curves over $\Q$. This requires some work; finally we obtain a parametrization of a subfamily of hyperelliptic curves: the parameters live on an elliptic curve over a function field. \\
We have implemented our algorithm by using many functions of GMP-ECM a free ECM program. During the computation, there are many multiplications by parameters of the hyperelliptic curves so the use of small parameters (i.e. which fit into one machine word) makes these multiplications negligible before the cost of full-length modular multiplications. This makes our software faster than GMP-ECM.

\bigskip

In Section $1$ we give some facts about ECM, decomposable genus $2$ hyperelliptic curves and Kummer surfaces. Details of the parametrization are given in Section $2$ while Section $3$ focuses on the HECM algorithm. In Section $4$ we describe our implementation and give some numerical results. 

\bigskip

In the following, $n$ denotes the number to be factored.

\section{background}

\subsection{ECM}

Following a common notation abuse in factoring algorithms, we work over $\Z/n\Z$ as if it were a field. The only operation that might fail is ``field'' inversion which is calculated using the Euclidean algorithm. If an inversion fails, we find a factor of $n$.

\bigskip

The ECM method starts by choosing a random elliptic curve $\mE$ over $\Z/n\Z$ and a point $P$ on it. Let $k$ be a positive integer, we compute $Q=[k]P$ and we hope $Q$ to be the zero of the curve modulo a prime $p$ dividing $n$ but not modulo $n$. If this is the case, then one division fails during the computation and we find the factor $p$. This happens if the cardinality $|\mE\left(\Fp\right)|$ divides $k$. There are two different phases (also called stages or steps) in ECM: in phase $1$ we hope that the cardinality is $B_{1}$-smooth (i.e. multiple of prime powers less than $B_{1}$) for some bound $B_{1}$. In phase $2$ we try to cover the ``close-miss'' case where the cardinality is a product of a $B_{1}$-smooth part by a prime cofactor. If this fails, we take another elliptic curve. \\
ECM is a probabilistic algorithm whose complexity is $ O(L(p)^{\sqrt{2}+o(1)}M( \log(n)) ) $ with $L(p)=e^{\sqrt{\log(p)\log(\log(p))}}$ and $M(\log(n))$ is the complexity of multiplications modulo $n$. The complexity of ECM is dominated by the size of the smallest factor $p$ of $n$ rather than the size of the number $n$ to be factored. Note however that ECM does not always find the smallest factor.

\bigskip

In stage $1$ we hope that the cardinality of the elliptic curve is $B_{1}$-smooth. We take $k=\prod_{ \pi\leq B_{1} } \pi^{ \left[  \log\left(B_{1}\right) / \log\left(\pi\right)  \right] }=\lcm\left(1,2,\ldots,B_{1}\right)$ so that all $B_{1}$-smooth numbers divide $k$. The main cost of stage $1$ is the arithmetic; different methods are used to reduce this cost. We focus on ECM with Montgomery formul\ae{} \cite{Montgomery87} since it shares several traits with the HECM method which is the subject of this work. Let $\mE$ be an elliptic curve in Montgomery's form: $y^2=x^3+a\,x^2+x$. We use projective coordinates $(x:y:z)$ to avoid divisions which are expensive. Moreover we disregard the $y$-coordinate which means that we identify a point with its opposite. It is still possible to ``double'' a point (i.e. computing $\pm2P$) and if the difference $\pm\left(P-Q\right)$ between two points is known then we can compute their sum $\pm\left(P+Q\right)$: this operation is called pseudo-addition. To compute the multiple of a point we need to find chains of doubling and pseudo-additions which is a special case of addition chains called ``Lucas Chains''~\cite{Montgomery83}. For instance
\[ 1\rightarrow2\rightarrow3\rightarrow4\rightarrow7\rightarrow10\rightarrow17\]
is a chain for $17$. One way to find such chains is to note that if we know $[n]P$ and $\left[n+1\right]P$ then we can compute $\left[2n\right]P$, $\left[2n+1\right]P$ and $\left[2n+2\right]P$. Hence we have binary chains: at each step we choose the point to double according to the binary expansion of the multiplier. For instance the following chain is the binary chain for~$17$: \[1\rightarrow2\rightarrow3\rightarrow4\rightarrow5\rightarrow8\rightarrow9\rightarrow17.\]
Binary chains are not the shortest chain of doubling and pseudo-addition. Montgomery's PRAC algorithm finds short Lucas chains \cite{Montgomery83}.

\bigskip

A description of stage $2$ falls out of scope of the present article since for reasons which will be explained later, stage $1$ is the only phase done by HECM.

\subsection{Decomposable curves}

\begin{definition}
A genus $2$ hyperelliptic curve $C$ is said to be decomposable if there exists an isogeny $\phi$ from $Jac\left(\mC\right)$ to the product of two elliptic curves $\mE_{i}$:
\[ \phi : \quad Jac\left(\mC\right) \longrightarrow \mE_{1} \times \mE_{2}. \]
If the kernel of the isogeny is $\Z/k\Z\times \Z/k\Z$ we say that the curve is $(k,k)$-decomposable.
\end{definition}

There are characterizations for a curve to be $(k,k)$-decomposable \cite{Shaska05}. We will only look at $\left(2,2\right)$-decomposable curves since the conditions are the simplest. 

\begin{theorem}
 Let $\mC$ be a hyperelliptic curve of genus $2$ given by the equation
\[ y^2=x\left(x-1\right)\left(x-\la\right)\left(x-\mu\right)\left(x-\nu\right) \]
where $\la$, $\mu$, $\nu$ are the Rosenhain invariants of the curve. $\mC$ is $\left(2,2\right)$-decomposable if and only if its Rosenhain invariants are linked by
\[ \la = \mu\frac{1-\nu}{1-\mu}. \]
\end{theorem}

The proof is given in \cite{GaSc01}. The two underlying elliptic curves have equation
\[ \chi\, y^2=\left(x-1\right)\left(x-x_2^2\right)\left(x-x_3^2\right)  \]
with
\[ q= \pm \sqrt{\mu\left(\mu-\nu\right)},\quad x_2= \frac{\mu+ q}{\mu- q}, \ x_{3}=\frac{1-\mu - q}{1-\mu + q}, \ \chi= -q\, \mu \left(\mu-1\right). \]
Over non algebraically closed fields $K$, the Rosenhain invariants may be non rational. Thus a $\left(2,2\right)$-decomposable curve can have a different equation. If $\mC$ is a $\left(2,2\right)$-decomposable curve in Rosenhain form then the two underlying elliptic curves are rational if and only if $\mu\left(\mu-\nu\right)$ is a square in $K$.\\
The morphisms between the hyperelliptic curves and the elliptic curves are defined on $K\left(q\right)$. They are given by
\[ \left(x,y\right) \longrightarrow \left( \left(\frac{ x-\mu-q}{x-\mu+q}\right)^2 , \frac{w \, y}{ \left(x-\mu+q\right)^{3}}  \right)  \ \text{with} \ w=\frac{8q}{\left(\mu-q\right)\left(-1+\mu-q\right)}. \]

\subsection{Kummer surfaces} \label{KS_gen}

Gaudry proposes in \cite{Gaudry07} to use theta functions to perform the arithmetic in the Jacobian of some genus $2$ curves. Working with theta functions allows the use of the numerous classical formul\ae{} that were found using complex analysis: the reader is referred to Mumford's books \cite{Mumford83,Mumford84} for details. In fact all the formul\ae{} we need are algebraic and can be used on finite fields (of characteristic different from $2$). \\
Let $\Omega$ be a matrix in the Siegel half space of dimension $2$. The theta functions are functions from $\C^{2}/\left(\Z^2+\Omega\Z^2\right)$ to $\C$. There are $16$ theta functions with half-integer characteristic: $10$ are even and $6$ are odd. The scalars obtained by evaluating them in $z=\left(0,0\right)$ are called theta constants. The main difference with Gaudry's article is that we use the squares of the coordinates. We are principally interested in only four squares of theta functions and the corresponding square theta constants $\abcd$.\\
 A Kummer surface $\K_{\abcd}$ can be defined by the choice of four ``theta constants'' $\abcd$ in $\PP^{3}\left(\C\right)$. We write $\XYZT$ for the coordinates of a point on $\K_{\abcd}$. The equation of the Kummer surface is: 
\begin{eqnarray*}
  4E'^{2} \al\be\ga\de XYZT &=& \left(\left(X^2+Y^2+Z^2+T^2\right)- \right. \\
 && \qquad \left. -F\left(X T+Y Z\right)-G\left(XZ+YT\right)-H\left(XY+ZT\right)\right)^2
\end{eqnarray*}
where $E'$, $F$, $G$, $H$ are constants that can be calculated from $\abcd$ by the formul\ae{} given in the appendix. 

\bigskip

Let $\frac{\epsilon}{\phi}$ be the ratio of two theta constants which can be calculated from the theta constants $\abcd$. Let
\[ \la:= \frac{\al \ga}{ \be \de}, \qquad \mu:= \frac{\ga \ep}{ \de \phi} , \qquad \nu:= \frac{\al \ep}{ \be \phi}, \]
\[ \mC:\qquad y^2=x\left(x-1\right)\left(x-\la\right)\left(x-\mu\right)\left(x-\nu\right). \]
Then $Jac\left(\mC\right)/\left\{\pm1\right\}$ is isomorphic to $\K_{\abcd}$ over $\C$:
\[ \psi: \qquad \J/\left\{\pm1\right\} \longrightarrow K_{\abcd} \]
Over non algebraically closed fields, let $\mCt$ be the quadratic twist of $\mC$. Then $\J/\left\{\pm1\right\}$ and $Jac(\tilde{\mathcal{C}})/\left\{\pm1\right\}$ are included in $\K_{\abcd}$:
\[ \K_{\abcd}= \psi\left(\J/\left\{\pm1\right\}\right) \bigcup \tilde{\psi}\left(\Jt/\left\{\pm1\right\}\right)\]
with the $2$-torsion points shared by the two Jacobians. Usually, we use the Mumford coordinates $\left(u,v\right)$ on $Jac\left(\mC\right)$, from a point $P$ in $K_{\abcd}$, it is possible to calculate $u$ and $v^2$ of its image by the morphism \cite{Gaudry07,Wamelen98}.
Of course $v$ is defined up to its sign since we can't distinguish between one point and its opposite. This generalizes the idea of Montgomery's formul\ae{} for elliptic curves. \\
The arithmetic on $Jac\left(C\right)$ transports to arithmetic on $\K_{\abcd}$. Given a point $P=\psi\left(D\right)$ on the Kummer surface it is possible to double it (i.e. to compute $\psi\left([2]D\right)$): see algorithm \ref{doubling}. Given two points on a Kummer surface $P=\psi\left(D\right)$ and $Q=\psi\left(D'\right)$, we can't add them since we don't know whether we must compute $\psi\left(D+D'\right)$ or $\psi\left(D-D'\right)$. However if one of these quantities is known then it is possible to compute the other, this is called a pseudo-addition (algorithm \ref{pseudo-addition}). We note $[2]P$ and $P+Q$ these two operations. Note that these formul\ae{} do not hold when one coordinate is zero.

\begin{algorithm}[t]
  \caption{Doubling on a Kummer surface}
  \label{doubling}
\begin{algorithmic}
	\REQUIRE a point $P=\XYZT$ on $\mathcal{K}_{\abcd}$.
	\ENSURE the point $[2]P=\left(X_{2}:Y_{2}:Z_{2}:T_{2}\right) $.
\begin{center}
\begin{tabular}{lllllll}
 $X'$ & $=$ & $\left(X+Y+Z+T\right)^2 \frac{1}{A}, $ & $\qquad$ & $Y'$ & $=$ & $\left(X+Y-Z-T\right)^2\frac{1}{B}, $ \\
 $Z'$ & $=$ & $\left(X-Y+Z-T\right)^2 \frac{1}{C}, $ & $\qquad$ & $T'$ & $=$ & $\left(X-Y-Z+T\right)^2\frac{1}{D} $ \\
 $X_{2}$ & $=$ & $\left(X'+Y'+Z'+T'\right)^2 \frac{1}{\al}, $ & $\qquad$ & $Y_{2}$ & $=$ & $\left(X'+Y'-Z'-T'\right)^2\frac{1}{\be}, $ \\
 $Z_{2}$ & $=$ & $\left(X'-Y'+Z'-T'\right)^2 \frac{1}{\ga}, $ & $\qquad$ & $T_{2}$ & $=$ & $\left(X'-Y'-Z'+T'\right)^2\frac{1}{\de} $ 
\end{tabular}
\end{center}
\end{algorithmic}
\end{algorithm}

\begin{algorithm}[t]
  \caption{Pseudo-addition on a Kummer surface}
  \label{pseudo-addition}
\begin{algorithmic}
	\REQUIRE two points $P=\XYZT$ and $Q=\left(\underline{X}:\underline{Y}:\underline{Z}:\underline{T}\right)$ on $\mathcal{K}_{\abcd}$, and the point $R=\left(\bar{X}:\bar{Y}:\bar{Z}:\bar{T}\right)$ equal to $P-Q$ such that $\bar{X}\bar{Y}\bar{Z}\bar{T} \neq 0$.
	\ENSURE the point $P+Q= \left(x:y:z:t\right) $.
\begin{center}
\begin{tabular}{lll}
 $X'$ & $=$ & $\left(X+Y+Z+T\right) \left(\underline{X}+\underline{Y}+\underline{Z}+\underline{T}\right) \frac{1}{A}, $ \\
 $Y'$ & $=$ & $\left(X+Y-Z-T\right) \left(\underline{X}+\underline{Y}-\underline{Z}-\underline{T}\right) \frac{1}{B}, $ \\
 $Z'$ & $=$ & $\left(X-Y+Z-T\right) \left(\underline{X}-\underline{Y}+\underline{Z}-\underline{T}\right) \frac{1}{C}, $ \\
 $T'$ & $=$ & $\left(X-Y-Z+T\right) \left(\underline{X}-\underline{Y}-\underline{Z}+\underline{T}\right) \frac{1}{D} $ \\
 $x$ & $=$ & $\left(X'+Y'+Z'+T'\right)^2 \frac{1}{\bar{X}}, $ \\
 $y$ & $=$ & $\left(X'+Y'-Z'-T'\right)^2 \frac{1}{\bar{Y}}, $ \\
 $z$ & $=$ & $\left(X'-Y'+Z'-T'\right)^2 \frac{1}{\bar{Z}}, $ \\
 $t$ & $=$ & $\left(X'-Y'-Z'+T'\right)^2 \frac{1}{\bar{T}} $
\end{tabular}
\end{center}
\end{algorithmic}
\end{algorithm}

We work with projective coordinates so divisions can be replaced by multiplications. Moreover the constants $1/\al$, $1/\be$, $1/\ga$, $1/\de$, $1/A$, $1/B$, $1/C$ and $1/D$ can be precomputed. The cost of pseudo-addition is $4$ divisions, $4$ multiplications, $4$ squares and $4$ multiplications by constants (thereafter denoted $4I+4M+4S+4d$) or if the divisions are replaced by multiplications the cost becomes $14M+4S+4d$. The cost of doubling is $8S+8d$. Note that by using the properties of projective coordinates, it is possible to save some multiplications \cite{Gaudry07}.

\bigskip

For the multiplication algorithm we want to avoid divisions since they are costly. This imposes that the inverses of the coordinates of $P-Q$ are always known when we want to pseudo-add $P$ and $Q$. Thus it is impossible to use ``short'' Lucas chains and the PRAC algorithm, we must use binary chains where $P-Q$ is always the initial point (algorithm \ref{Multiplication}). For each bit of the multiplier $k$, we do one doubling and one pseudo-addition on the Kummer surface. The total cost of a multiplication by~$k$ is $7M+9S+9d$ per bit of $k$ since some computations can be shared between the two operations (it is $12M+4S+16d$ if we don't share them).

\begin{algorithm}[t]
  \caption{Multiplication algorithm}
  \label{Multiplication}
\begin{algorithmic}
	\REQUIRE a point $P$ on $\mathcal{K}_{\abcd}$ and an integer $k$. 
	\ENSURE the point $[k]P $.
	
	\IF{$k=2$} \RETURN $[2]P$
	\ELSE		
		\STATE Let $k=\sum_{i=0}^{l}k_{i}2^{l-i}$ be the binary expansion of $k$ with $k_{0}=1$ the most significant bit.
		\STATE $P_{m}=P$; $P_{p}=[2]P$;

		\FOR{$i$ from $2$ to $l$}

			\STATE $Q=P_{p}+P_{m}$  \COMMENT{note that we have $P_{p}-P_{m}=P$}

			\IF{$k_{i}=1$}
				\STATE $P_{p}=[2]P_{p}$; $P_{m}=Q$;
			\ELSE[$k_{i}=0$]
				\STATE $P_{m}=[2]P_{m}$; $P_{p}=Q$;
			\ENDIF

		\ENDFOR

		\RETURN $P_{m}$

	\ENDIF
\end{algorithmic}
\end{algorithm}

\begin{remark}\label{remark-coor_ini_point}
The coordinates of the initial point $P$ affect the cost of computation. For instance if two coordinates of $P$ are equal or opposite then we gain, for each bit of the multiplier, one multiplication and another one which in fact is a square. This is a speed up of $5\%$.
\end{remark}

\section{Parametrization} \label{para}

In this Section, we find suitable curves for HECM. We obtain curves over $\Q$ having the desired properties and consider their reduction modulo $n$. We want to have an infinite number of curves over $\Q$ so that we have enough curves over $\Z/n\Z$. \\
Note that throughout the construction, we should never compute a square root because of two reasons: first this square root could be in an extension of $\Q$, thus when reducing modulo $n$ we could be working in an extension of $\Z/n\Z$ and the arithmetic would be slower. Moreover, computing square roots modulo the composite integer $n$ is equivalent to find factors of $n$. To sum up, we require that all constants are defined by rational functions in several variables.

\subsection{Parametrization of the hyperelliptic curves}

\begin{lemma}
Let $\mC$ be a genus $2$ curve over $\Q$ with equation
\[ \mC :\qquad \chi y^2 = x\left(x-1\right)\left(x-\la\right)\left(x-\mu\right)\left(x-\nu\right). \]
$\mC$ can be used in HECM (i.e. $\mC$ is $(2,2)$-decomposable with rational underlying elliptic curves and the rational Kummer surface) if and only if
\[ \la=\mu\frac{1-\nu}{1-\mu}, \qquad \mu\left(\mu-\nu\right) = \square, \qquad  \la \mu \nu = \square \]
where $\square$ means that the quantity must be a square.
\end{lemma}

\begin{proof}
The first two conditions mean that $\mC$ is $(2,2)$-decomposable with rational underlying elliptic curves. The Kummer surface is rational if and only if the squares $\abcd$ of the four theta constants are rational. They are linked with the Rosenhain invariants by:
\[ \la= \frac{\al \ga}{ \be \de}, \qquad \mu= \frac{\ga \ep}{ \de \phi} , \qquad \nu= \frac{\al \ep}{ \be \phi} \]
where $\ep$ and $\phi$ are two other squares of theta constants whose ratio is rational when $\abcd$ are rational. Note that this implies that the Rosenhain invariants must also be rational so it gives a justification to the choice of the equation of the hyperelliptic curve. From these three equations, we obtain:
\[ \frac{\al}{\be}=\frac{\sqrt{\la\mu\nu}}{\mu}, \qquad \frac{\ga}{\de}=\frac{\sqrt{\la\mu\nu}}{\nu}. \]
So the product $\la\mu\nu$ must be a square in $\Q$. \\
When we write the equation $\la=\mu\frac{1-\nu}{1-\mu}$ in terms of theta constants (use the formul\ae{} from \cite{GaSc01} and \cite{Wamelen98}), it yields to $\al^2=\de^2$. Moreover $\be$, $\ga$, $\de$ are functions of $\la$, $\mu$, $\nu$ and $\al$ and thus are rational. This proves that the conditions are sufficient. 
\end{proof}

\begin{remark}
The Rosenhain invariants $\left(\la,\mu,\nu\right)$ are defined up to the action of the group $PGL(2,5)$ see \cite{GaSc01}. Our choice here is one that leads to an equality between two of the first four theta constants. With a different choice of order, we would have had another square root to handle. Moreover the fact that $\al^{2}=\de^{2}$ is a main advantage for the arithmetic: it saves two scalar multiplications per multiplier's bit when computing the multiple of a point.
\end{remark}

Since $\abcd$ live in $\PP^{3}\left(\Q\right)$, we can take $\al=1$. We choose also $\de=\al=1$. The choice $\de=-\al=-1$ leads to an isomorphic Kummer surface. \\
Take $ \mu=1-\frac{\nu\left(1-\nu\right)}{s^2}$ with $s\in\Q$ so that $ \la\mu\nu= \mu^2 s^2$. The second equation becomes
\[ \mu\left(\mu-\nu\right) = \frac{1}{s^4} \left(-\nu+\nu^2+s^2\right)\left(\nu-1\right)\left(\nu-s^2\right) = \square. \]
Assume that $\frac{\nu-s^2}{\nu-1}$ is a square $u^2$ (then $\nu=\frac{s^2-u^2}{1-u^2}$) so that the equation above rewrites as:
\[ 1 + (-3/s^2+1/s^4)u^2+u^4/s^2=\square. \]
Let $v^2$ be the square, the point $\left(u,v\right)$ lies on an elliptic curve over $\Q\left(s\right)$. This elliptic curve is in the Jacobi model \cite{Duquesne07}, is of rank $1$ and we have a non torsion point $P=\left(1,1-\frac{1}{s^2}\right)$ on it. 

\begin{theorem}
A subfamily of $(2,2)$-decomposable hyperelliptic genus $2$ curves with rational underlying elliptic curves and with rational Kummer surfaces is given by the following parametrization:
\[ \mC :\qquad \chi y^2 = x\left(x-1\right)\left(x-\la\right)\left(x-\mu\right)\left(x-\nu\right), \]
\[ \la=\mu\frac{1-\nu}{1-\mu},\qquad \mu=1-\frac{\nu\left(1-\nu\right)}{s^2}, \qquad \nu=\frac{s^2-u^2}{1-u^2} \]
where $\left(u,v\right)$ lies on the following elliptic curve
\[ 1 + (-3/s^2+1/s^4)u^2+u^4/s^2=v^2. \]
A non-torsion point on this curve is $(1,1-\frac{1}{s^2})$. The parameters $s$, $\chi$, $u$, $v$ are rational numbers that must verify the following conditions.
\end{theorem}

\begin{condition}\label{cond_hcurve}
With the notations of the theorem, the curve $\mC$ is of genus $2$ if and only if $0$, $1$, $\infty$, $\la$, $\mu$ and $\nu$ are distinct and $\chi$ is not zero. This is equivalent to 
\[ \chi \neq 0, \qquad s\neq0,\pm1, \qquad u \neq0,\pm1, \qquad v\neq0, \]
\[ s\neq\pm u, \qquad s^2-2u^2+u^4\neq0. \]
\end{condition}

In particular, condition \ref{cond_hcurve} implies that we can't use the point $(1,1-\frac{1}{s^2})$ on the Jacobi curve. We must take a multiple of this point: for instance its double $(2,1+\frac{2}{s^2})$.

\begin{condition}\label{cond_Kummer} 
For the arithmetic on the Kummer surface, all the theta constants must be non zero. The parameters must satisfy condition \ref{cond_hcurve} and 
\[ s\neq \pm u^2 . \]
\end{condition}

\bigskip

The parameter $\chi$ determines whether we are on the curve or on its quadratic twist. However it is not chosen during the parametrization of the curve: since the curve and its twist are both included in the same Kummer surface, the choice of a point on it determines whether we work on the curve or on its twist. In practice, the $y$-coordinate of the point is never used, so $\chi$ is never computed. 

\bigskip

We could have taken $\frac{\nu-s^2}{\nu-1}=l\,u^2$ (here $l=1$) but that would have led to more complicated equations. For instance with $l=-1$, we found an elliptic curve with no point on it. In this case ($l=-1$), we had to assume that $s$ was on a conic for the curve to become of rank $1$.

\begin{remark}
ECM supposes that the cardinality of elliptic curves on $\Fp$ behaves like a random integer of size around $p$ with some additional divisibility conditions. In HECM, we work simultaneously with two elliptic curves. Their cardinalities are not independent: points of $4$-torsion cannot exist on only one curve. However, experimentally, it seems that this is the only noticeable link between their cardinalities. 
\end{remark}

\subsection{Finding a point on the Kummer surface}

In ECM we need a generic point (i.e. a non torsion point) in the group. Contrary to the Brent-Suyama parametrization for elliptic curves, it is possible here to find a point after having chosen a hyperelliptic curve. However, the initial point $P$ must be a non torsion point of the two underlying elliptic curves over $\Q$ and not only of the Jacobian of the hyperelliptic curve. Remember that we don't want to take square roots or to work in an extension field so we can't take three random coordinates and solve the equation of the Kummer surface to find the last one.

\bigskip

The first method to find a point is to take one coordinate to be zero. Then the three other lie on a conic. Since conics are birationally equivalent to $\PP^{1}$, we have an infinite number of points on the Kummer surface. In general, these points are not torsion points. Note however that these points cannot be used directly for the multiplication algorithm: we must double them until we find a point with no zero coordinate; in general, one doubling is enough. \\ 
There is a better solution to find points: we look for points with two equal or opposite coordinates. Such points have the advantage of saving multiplications (see remark \ref{remark-coor_ini_point}): the cost of the multiplication algorithm becomes $12M+10S$ per bit of the multiplier. It turns out that not all choices of coordinate pairs are possible: some choices lead only to torsion points or to impossible equations like $-1=\square$. The choice $Y=-X$ yields to the equation of an elliptic curve with many points on it. These points lead to points on the Kummer surface which are of infinite order (over $\Q$). For instance we use the following point: 
\begin{eqnarray*}
X =&  \frac{u^2\left(u^2-1\right)\left(su^{4}-3su^{2}+u^{2}+s^{3}-s^{2}+s\right)}{\left(s-u^2\right)s^{4}v^{2}}, \qquad &Y = -X \\
T =&  \frac{\left(s^2-u^2\right)\left(u^2-1\right)}{s^{4}v^{2}}, \qquad & Z = 1.
\end{eqnarray*}

\section{Algorithm and ameliorations}

\subsection{HECM}

In ECM, the big product $k=\prod_{ \pi\leq B_{1} } \pi^{ \left[  \log\left(B_{1}\right) / \log\left(\pi\right)  \right] }$ is not computed as such \cite{ZiDo06}. Instead, for each prime $\pi$ we compute $l=\left[\log\left( B_{1}\right) / \log\left(\pi\right)\right]$ and then we do $l$ times $Q = [\pi]Q$ and go to the next prime. This avoids the computation of a big integer product and PRAC finds shorter chains if it works with prime multipliers~\cite{Montgomery83}. If we work with Kummer surfaces we can't use PRAC because changing the initial point is costly (we would have divisions). Moreover the advantage of choosing a good initial point is lost if we don't use binary Lucas chains. Therefore in HECM we begin by computing $k=\prod_{ \pi\leq B_{1} } \pi^{ \left[  \log\left(B_{1}\right) / \log\left(\pi\right)  \right] }$. In theory this computation is costly but it practice it is negligible: we use product tree and fast multiplication. Moreover this computation is shared for many curves. 

\bigskip

We hope to encounter the zero of one of the underlying elliptic curves so it is important to have explicit morphisms between the Kummer surface and the two underlying elliptic curves. In fact, it is impossible to obtain real points on the curves since we can't distinguish between a point and its opposite and we don't even know if we are on the curves or on their twists. If the elliptic curves are in the Weierstrass form, we only need the coordinates $\left(x::z\right)$ to test if the point is zero: just compute $\gcd\left(z,n\right)$. The morphisms between the Kummer surface and the underlying elliptic curves (identified with their quadratic twists) are rational over $\Q$. However their building blocks are defined over an extension field of $\Q$: they use square roots that disappear in the global morphisms. Since we are not allowed to use square roots we must rewrite the global morphisms. As a result their expressions become complicated however, they are computed only once for every run of HECM, thus they are negligible before the cost of computing $[k]P$. Note that we only need the $(x::z)$ coordinates which are invariant under the isomorphism from the curve (in the Weierstrass form) to its twist. See the appendix for a description of the different elementary blocks. The first stage of the HECM algorithm is summarized in algorithm \ref{HECM}.

\bigskip

We now turn to what is called stage $2$. The initial point for the arithmetic is $Q=[k]P$, thus we don't have the benefit of a ``good'' initial point. Moreover, stage $2$ needs the $x$-coordinate of many points on the elliptic curve in Weierstrass form~\cite{ZiDo06}, therefore if we use hyperelliptic curves we need to apply morphisms a lot and, in this case, their cost would not be negligible. In addition to that, if we want to use Brent-Suyama's extension, we can't use pseudo-additions but real additions. For all these reasons, it seems that hyperelliptic curves should not be used for stage~$2$. Instead, we can apply the ECM stage $2$ to the two underlying elliptic curves.

\begin{algorithm}[t]
  \caption{HECM (stage 1)}
  \label{HECM}
\begin{algorithmic}[1]
	\REQUIRE the number $n$ to be factor. The smoothness bound $B_{1}$.
	\ENSURE a factor $p$ of $n$.
	
	\STATE Compute $ k=\lcm\left(1,2,\ldots,B_{1}\right) $.
	\STATE Choose a random decomposable curve $\mC$ over $\Z/n\Z$ and a point $P$ on its Kummer surface.
	\STATE Compute $Q=[k]P$.
	\STATE Map $Q$ to the two underlying elliptic curves $\mE_{i}$.
	\STATE Hope that $Q=\mathcal{O}_{\left(\mE_{i}\right)} \mod p$ for one $\mE_{i}$ (test whether $\gcd\left(z,n\right)\neq1$).
	\STATE Else go to $2$.

\end{algorithmic}
\end{algorithm}

\subsection{Torsion of the curves} \label{torsion}

To improve the probability of success, we can force the group order of the elliptic curves modulo primes $p$ to be divisible by small numbers. Mazur's theorem states that the torsion group $E_{tor}\left(\Q\right)$ of any elliptic curve over $\Q$ is isomorphic to one of the following groups:
\[ E_{tor}\left(\Q\right)\cong \left\{ 
\begin{array}{lll}
 \Z/m\Z                & \qquad & 1\leq m\leq 10\ or\ m=12 \\
 \Z/2\Z \times \Z/2m\Z & \qquad & 1\leq m\leq4
\end{array}
\right. \]
In our case the underlying elliptic curves have at least four $2$-torsion points which means that $\Z/2\Z \times \Z/2\Z$ is a subgroup of the torsion group. 

\bigskip

In theory it should be possible to find decomposable hyperelliptic curves with underlying elliptic curves having larger torsion group. However this leads to more complicated equations for the parametrization. There is also a theoretical obstruction to improving the torsion: to have a torsion point of order more than $2$ means to be able to find a starting point on the Kummer surface corresponding to a point on $\J$ and not on its twist $Jac(\tilde{\mathcal{C}})$. This theoretical obstruction could be removed if the curve and its twist had the same torsion. However this increases the number of equations needed for the parametrization.

\bigskip

Using the parametrization presented in Section \ref{para}, we have table \ref{4_torsion} where letters indicate whether and where $4$-torsion points exist: $C$ for the curve and $T$ for its quadratic twist. With the same probability, we work on the curve or on its twist. Suppose that the Legendre symbols of $s^2-u^2$, $s^2-1$, $u^2-1$ are independent which experimentally is a valid approximation. Then if $p\equiv1\, \left[4\right]$ we have a point of $4$-torsion with probability $1/2$, and if $p\equiv3 \, \left[4\right]$, the probability is $3/4$. Computer experiments show that for $p\equiv3 \, \left[4\right]$ the power of two dividing the cardinality is, on average, $3.48$ (instead of $3.5$) and for $p\equiv1\, \left[4\right]$ it is $3.15$ (instead of $3$).

\bigskip

\begin{table}
\caption{$4$-torsion points on the underlying elliptic curves}
\label{4_torsion}
\begin{center}
 \begin{tabular}{|c|c|c|cc|cc|}
\hline
                                        &                                         &                                        & \multicolumn{4}{|c|}{$-1$} \\
\cline{4-7}
\raisebox{1.5ex}[0cm][0cm]{$s^2-u^2$}     & \raisebox{1.5ex}[0cm][0cm]{$s^2-1$}     &\raisebox{1.5ex}[0cm][0cm]{$u^2-1$} & \multicolumn{2}{|c}{$\square$} & \multicolumn{2}{|c|}{$\boxtimes$}   \\
\hline

                                        &                                         & $\square$                              & $C$ &     & $C$ & $T$ \\
\cline{3-7}
                                        & \raisebox{1.5ex}[0cm][0cm]{$\square$}   & $\boxtimes$                            & $C$ &     & $C$ & $T$ \\
\cline{2-7}
\raisebox{1.5ex}[0cm][0cm]{$\square$}   &                                         & $\square$                              & $C$ &     &     &  					    \\
\cline{3-7}
                                        & \raisebox{1.5ex}[0cm][0cm]{$\boxtimes$} & $\boxtimes$                            &     & $T$ & $C$ & $T$ \\
\hline
                                        &                                         & $\square$                              &     & $T$ & $C$ & $T$ \\
\cline{3-7}
                                        & \raisebox{1.5ex}[0cm][0cm]{$\square$}   & $\boxtimes$                            & $C$ &     &     & 
					    \\
\cline{2-7}
\raisebox{1.5ex}[0cm][0cm]{$\boxtimes$} &                                         & $\square$                              & $C$ &     & $C$ & 						$T$ \\
\cline{3-7}
                                        & \raisebox{1.5ex}[0cm][0cm]{$\boxtimes$} & $\boxtimes$                            & $C$ &     & $C$ &   $T$ \\
\hline
 \end{tabular}
\end{center}
\end{table}

Experimentally, the order of the curves is as likely to be smooth as a random integer about $1/7.75$ in value (compared with $1/23.7$ for Suyama's curves \cite{ZiDo06}). \\
The parametrization $s=\frac{3+t^2}{3-t^2}$, $u=2$, $v=1+\frac{2}{s^2}$ provides better torsion (the power of two is on average $3.66$) but there are few such curves with low parameters.

\subsection{Small parameters}

For a chosen Kummer surface and a chosen initial point, the arithmetic on Kummer surfaces uses many multiplications by fixed parameters. Remember from section \ref{KS_gen} that the cost of a scalar multiplication is $4M+12S+16d$ per bit of the multiplier; i.e., $4$ multiplications, $12$ squares and $16$ multiplications by constants. Suppose that these constants are small in comparison to the number $n$ to be factored, then the cost of multiplications by them will be negligible with respect to the cost of the full-length multiplications. The cost per bit of the multiplier becomes $4M+12S$ and if we use the rule of thumb $1S=0.8M$ then it is $13.6M$. For one run, GMP-ECM (which uses Montgomery coordinates and the PRAC algorithm \cite{Montgomery83}) uses approximately $6M+3S$. Thus two runs of ECM cost $12M+6S$ or $16.8M$ with $S=0.8M$. This is a speed-up of $20\%$. Of course this is only valid for large $n$. \\ 
EECM (ECM with Edwards curves) uses signed sliding window to compute $\left[k\right]P$ (see \cite{BeBiLaPe08}). These chains use $1$ elliptic curve doubling and $\epsilon$ elliptic curve additions for each bit of $k$ where $\epsilon$ converges to $0$ when $k$ increases. Doubling in Edwards coordinates uses only $3M+4S$, while additions use $10M+1S+1d$. In theory, two runs of EECM would be faster than one run of HECM if $\epsilon$ is less than $1/12$ (if $1S=0.8M$ then $\epsilon$ should be less than $1/20$). However this requires that the width of the window be very large and thus precomputations and memory usage might be not negligible. The authors of \cite{BeBiLaPe08} claim that with $B_{1}=16384$, EECM uses $195111$ multiplications in stage $1$ (with a width-$6$ signed sliding window). For the same $B_{1}$, HECM uses $379334$ multiplications but since it does two curves in parallel, it uses only $189667$ multiplications for stage $1$ for one elliptic curve. Therefore, more experimentations are needed but to our knwolegde there is no public implementation of EECM. 

\bigskip

By a small constant we mean a number which fits into a signed long. This constraint limits the number of curves we can use: if we work with a 64 bits processor there are $185,399$ useful hyperelliptic curves which is sufficient for finding factors of more than 65 digits. In the multiplication algorithm, the $16d$ could be reduced to $12d$ since we work in the projective space but we would have rational numbers on $\Q$ that, modulo $n$, must each fit into a long which is more complicated. \\
More explicitly, our parametrization gives $\abcd$ in terms of rational functions of small degree in $(s,u)\in\Q^{2}$. The same is true for $\left(1/A:1/B:1/C:1/D\right)$. We work in projective coordinates so we can clear the denominators. The constants used in the multiplications become:
\[\left(\frac{1}{\al}:\frac{1}{\be}:\frac{1}{\ga}:\frac{1}{\de}\right)=\left(s^{4}v^{2}:\ s^{5}v^{2}:\ s\left(u^2-s^2\right)\left(u^2-1\right):\ s^{4}v^{2} \right), \]
\begin{eqnarray*}
\left(\frac{1}{A}:\frac{1}{B}:\frac{1}{C}:\frac{1}{D}\right) &=& \left(\left(s-1\right)^2\left(s^2-2u^2+u^4\right)\left(s+u^2\right)^2: \right.\\
& & \left. \left(u^2-1\right)\left(s+u^2\right)^2\left(s-u^2\right)^2: \right. \\
& & \left. -\left(u^2-1\right)\left(s+u^2\right)^2\left(s-u^2\right)^2:\right. \\ 
& & \left. \left(s+1\right)^2\left(s^2-2u^2+u^4\right)\left(s-u^2\right)^2 \right).
\end{eqnarray*}
Note that $s$ is rational so we take $s=a/b$ for integers $a$,$b$. After generating a hyperelliptic curve, we test if these constants are small. We then need to find a point on the Kummer surface for which the inverses of its coordinates are small. Unfortunately we have found no point such that the polynomial expressions for the inverses of the coefficients have degree equal or lower than the degree of the constants above. In practice, we have several generic points so we check whether those points fit, if not we try with another curve.

\subsection{A numerical example}

Let's try to factor $n=4,816,415,081$ with HECM. We use $B1=25$ and $B2=200$. Let's use the following parameters: $s=\frac{1}{2}$, $u=2$ and $v=9$. The Kummer surface is given by the parameters
\[ \abcd = \left(1:2:\frac{9}{10}:1\right) . \]
One point on the Kummer surface is 
\[ P=\left(-272:272:63:-140\right) . \]
We compute $\left[k\right]P$ on the Kummer surface with $k=\lcm\left(2,B1\right)=26,771,144,400$ and send the resulting point on the two underlying elliptic curves:
\[ \left(3455587574,1\right) \text{ on } 734346861\,y^2=\left(x-1\right)\left(x-\frac{1}{25}\right)\left(x-\frac{1}{121}\right), \]
\[ \left(3222355131,1\right) \text{ on } 2791313056\,y^2=\left(x-1\right)\left(x-25\right)\left(x-121\right). \]
We now use the ECM stage 2 on the two elliptic curves. The first one does not produce any factor but the second one produces the number $83003$ which is a factor of $n=4816415081=58027*83003$. \\
We inspected the calculus and found that we didn't work on the hyperelliptic curve but on its twist:
\[ \chi y^2=x\left(x-1\right)\left(x-\frac{9}{20}\right)\left(x-\frac{9}{4}\right)\left(x-\frac{5}{4}\right)\]
where $\chi$ is a non quadratic residue. The order of the initial point in the Jacobian of this curve modulo $83003$ is $2 \cdot 3\cdot5\cdot11 \cdot19\cdot73\cdot631$. Its order on the first elliptic curve is $2\cdot 7\cdot631$ and on the second $2 \cdot 3\cdot5 \cdot19\cdot73$.

\section{Implementations and results}

There are many implementations of ECM: GMP-ECM is a free program based on the GNU MP library. It is described in detail in \cite{ZiDo06}. By using many functions from GMP-ECM, we have built a new program called GMP-HECM which uses hyperelliptic curves. GMP-HECM will be distributed with GMP-ECM (currently it is distributed with the development version at http://gforge.inria.fr/projects/ecm/). Our software does stage $1$ by generating a genus $2$ hyperelliptic curve with small parameters, computes $\left[k\right]P$ on the Kummer surface and send the resulting point on the two underlying elliptic curves. Then, for stage $2$, it uses GMP-ECM for the two underlying elliptic curves (This could be done in parallel). Indeed, GMP-ECM can use the results of stage $1$ from an other software: it just needs the parameter $A$ of the elliptic curve in short Weierstrass form $y^2=x^3+Ax+B$ and the coordinates of the point on it. 

\bigskip

Our elliptic curves are not, in general, Suyama curves so they don't have the same torsion groups as it is expected in GMP-ECM. This affects the probability of the cardinality of the curves being smooth and thus the choice of $B_{1}$ and $B_{2}$. GMP-ECM is able to choose the parameters $B_{2}$ for stage $2$ from the value of $B_{1}$. This choice is not the same as the ``standard continuation'' (i.e. $B_{2}=100*B_{1}$) but produces much bigger $B_{2}$. This reduces the number of expected curves for finding a factor. The choice of $B_{1}$ is heuristic, we tried to minimize the means of the expected time. Note that for finding small factors, HECM is not interesting due to the cost of initialization and morphisms. Table \ref{choice_B1} compared the optimal value of $B_{1}$ for different sizes of the (usually unknown) factor. We see that for finding factors of at least $35$ digits, we use the same $B_{1}$ than ECM. Moreover the ratio of the expected number of curves needed for finding a factor decreases with the size of the factor.

\begin{table} 
\caption{Optimal choice of the parameters for ECM and HECM}
\label{choice_B1}
\begin{center}
\begin{tabular}{|c||r|c||r|r|c|}
\hline
Size    & \multicolumn{2}{|c||}{GMP-ECM}                    &   \multicolumn{3}{|c|}{GMP-HECM}           \\
\cline{2-6}
of the  & \multicolumn{1}{c|}{Optimal} & Expected number    & \multicolumn{1}{c|}{Optimal} & \multicolumn{1}{c|}{Default} & Expected number \\
factor  & \multicolumn{1}{c|}{$B_{1}$} & of elliptic curves & \multicolumn{1}{c|}{$B_{1}$} & \multicolumn{1}{c|}{$B_{2}$} & of elliptic curves \\ \hline
$10^{20}$   &                 $11,000$ & $74$               &                     $14,000$ &        $2.10^{6\phantom{0}}$ & $75$ \\
\hline
$10^{25}$   &                 $50,000$ & $214$              &                     $60,000$ &       $16.10^{6\phantom{0}}$ & $214$ \\
\hline
$10^{30}$   &                $250,000$ & $430$              &                    $260,000$ &      $130.10^{6\phantom{0}}$ & $491$ \\
\hline
$10^{35}$   &                 $1.10^6$ & $904$              &                    $1.10^6 $ &      $900.10^{6\phantom{0}}$ & $1,116 $ \\
\hline
$10^{40}$   &                 $3.10^6$ & $2,350$            &                     $3.10^6$ &        $4.10^{9\phantom{0}}$ & $2,871$ \\
\hline
$10^{45}$   &                $11.10^6$ & $4,480$            &                    $11.10^6$ &       $28.10^{9\phantom{0}}$ & $5,425$ \\
\hline
$10^{50}$   &                $43.10^6$ & $7,553$            &                    $43.10^6$ &      $200.10^{9\phantom{0}}$ & $9,003$ \\
\hline
$10^{55}$   &               $110.10^6$ & $17,769$           &                   $110.10^6$ &      $750.10^{9\phantom{0}}$ & $21,183$ \\
\hline
$10^{60}$   &               $260.10^6$ & $42,017$           &                   $260.10^6$ &                  $2.10^{12}$ & $49,534$ \\
\hline
$10^{65}$   &               $850.10^6$ & $69,408$           &                   $850.10^6$ &                 $14.10^{12}$ & $81,387$ \\
\hline
\end{tabular}
\end{center}
\end{table}

\bigskip

For not too small $B_{1}$, the theoretical and the real cost of the algorithm is linear in $B_{1}$. However its dependence in the size of the input number is not simple: modu\-lar multiplications are quadratic for small $n$ but become quasi-linear for large~$n$. Moreover, there are a lot of other operations which complexity is not negligible for small $n$. Table \ref{time_operations} presents the time taken by the different operations during the computation of $[k]P$ for different sizes of $n$. We see that squares take the same time as multiplications. The reason is that GMP-ECM has special assembly code for general modular multiplications but not for squaring and thus squarings use the same assembly code as normal multiplications. For multiplications by small constants, we made special assembly code. The table also shows that additions are not negligible for small $n$. 

\bigskip

Table \ref{comparison} compares the running time of two runs of GMP-ECM versus GMP-HECM (remember that one run of HECM is equivalent to two runs of ECM) for a fixed $B_{1}$ and different sizes of $n$. This shows that for large $n$ (at least $10^{250}$), our software GMP-HECM is faster than GMP-ECM. Since HECM uses more squarings than ECM, having optimized assembly squaring code would be likely to provide a small improvement to the limit where GMP-HECM is more interesting than GMP-ECM. For very large $n$ the speed up of GMP-HECM compared to GMP-ECM is around $11\%$. This result agrees with the theoretically predicted value.

\begin{table} 
\caption{Fraction of time taken by the arithmetic operations for different sizes of $n$}
\label{time_operations}
\begin{center}
\begin{tabular}{|c||c|c|c|c|c|c|c|c|c|}\hline
size of $n$ & $10^{100}$ & $10^{150}$ & $10^{200}$ & $10^{250}$ & $10^{300}$ & $10^{350}$ & $10^{400}$ & $10^{500}$ & $10^{1000}$ \\ \hline\hline
M           & $0.143$    & $0.159$    & $0.194$    & $0.216$    & $0.209$    & $0.219$    & $0.216$    & $0.222$    & $0.234$     \\ \hline
S           & $0.428$    & $0.465$    & $0.557$    & $0.597$    & $0.622$    & $0.645$    & $0.65$     & $0.682$    & $0.700$     \\ \hline
d           & $0.163$    & $0.127$    & $0.092$    & $0.081$    & $0.070$    & $0.057$    & $0.037$    & $0.031$    & $0.051$     \\ \hline
additions   & $0.237$    & $0.179$    & $0.116$    & $0.109$    & $0.074$    & $0.057$    & $0.059$    & $0.041$    & $0.017$     \\ \hline\hline
S/M         & $1.00$     & $0.98$     & $0.96$     & $0.93$     & $0.99$     & $0.98$     & $1.04$     & $1.02$     & $1.00$      \\ \hline
d/M         & $0.26$     & $0.20$     & $0.12$     & $0.093$    & $0.084$    & $0.066$    & $0.043$    & $0.035$    & $0.055$     \\ \hline
\end{tabular}
\end{center}
\end{table}

\begin{table} 
\caption{Comparison between stage $1$ of GMP-ECM and GMP-HECM for different sizes of $n$ with a constant $B_{1}=10^7$ on a core 2 at 2.4Ghz.}
\label{comparison}
\begin{center}
\begin{tabular}{|c||c|c|c|c|c|c|c|c|c|}\hline
size of $n$          & $10^{100}$ & $10^{150}$ & $10^{200}$ & $10^{250}$ & $10^{300}$ & $10^{350}$ & $10^{400}$ & $10^{500}$ & $10^{1000}$ \\ \hline
                     &            &            &            &            &            &            &            &            &              \\
$\frac{HECM}{2*ECM}$ & $1.15$     & $1.09$     & $1.03$     & $1.0$      & $0.97$     & $0.94$     & $0.95$     & $0.93$     & $0.89$       \\ 
                     &            &            &            &            &            &            &            &            &              \\
\hline 
\end{tabular}
\end{center}
\end{table} 

\section*{Conclusion}

Using a subfamily of hyperelliptic curves of genus $2$ (i.e decomposable curves), we have built an algorithm for factorization which is equivalent to two simultaneous runs of ECM. Kummer Surfaces allowed us to choose small parameters that give an efficient arithmetic in the Jacobian of the curve. In practice, for large numbers, our implementation is faster than GMP-ECM. Special assembly code for squaring would provide another improvement to GMP-HECM. 

\bigskip

Many thanks to Emmanuel Thom\'{e} for his help on this work and for his comments on the draft versions. Thanks also to the developers of GMP-ECM and in particular to Alexander Kruppa for his great assembly code.

\section*{Appendix: Morphisms}

A Kummer surface is defined by the first four theta constants $\abcd$. Its equation is
\begin{eqnarray*}
  4E'^{2} \al\be\ga\de XYZT &=& \left(\left(X^2+Y^2+Z^2+T^2\right)- \right. \\
 && \qquad \left. -F\left(X T+Y Z\right)-G\left(XZ+YT\right)-H\left(XY+ZT\right)\right)^2
\end{eqnarray*}
where
\begin{center}
\begin{tabular}{lllllll}
 $A$ & $=$ & $\al+\be+\ga+\de, $ & $\qquad$ & $Y'$ & $=$ & $\al+\be-\ga-\de, $ \\
 $C$ & $=$ & $\al-\be+\ga-\de, $ & $\qquad$ & $T'$ & $=$ & $\al-\be-\ga+\de, $ \\
\end{tabular}
\end{center}
\begin{eqnarray*}
E' & = & \frac{ABCD}{\left(\al\de-\be\ga\right)\left(\al\ga-\be\de\right)\left(\al\be-\ga\de\right)} , \\
F  & = & \frac{\left(\al^2-\be^2-\ga^2+\de^2\right)}{\left(\al\de-\be\ga\right)} , \\
G  & = & \frac{\left(\al^2-\be^2+\ga^2-\de^2\right)}{\left(\al\ga-\be\de\right)} , \\
H  & = & \frac{\left(\al^2+\be^2-\ga^2-\de^2\right)}{\left(\al\be-\ga\de\right)} .
\end{eqnarray*}

In our case (decomposable curves with $\al=\de$), the other theta constants are given by the following relations:
\[ \tun\zero=\al, \quad \tde\zero=\be, \quad \ttr\zero=\ga,\quad \tqu\zero=\de,\]
\[ \thu\zero=\ep=\sqrt{\al\be-\ga\de}\sqrt{\frac{\be}{\al}}\sqrt{\nu}=\sqrt{\al\be-\ga\de} \sqrt{\frac{\sqrt{\la\mu\nu}}{\la}},  \]
\[ \tdi\zero=\phi=\frac{\ep\al}{\nu\be}=\frac{\al\be-\ga\de}{\ep}, \]
\[ \tse\zero=\sqrt{-\phi^2}, \]
\[ \tne\zero=\sqrt{-\ep^2}, \]
\[ \tci\zero=\frac{\ga\ep-\de\phi}{\tne\zero}, \]
\[ \tsi\zero=\frac{\al\de-\be\ga}{\tci\zero}. \]

\bigskip

The global morphism in HECM goes from the Kummer surface to the product of the underlying elliptic curves modulo the isomorphisms which identify a curve with its twist. However we only need the $(x::z)$ coordinates which are invariant under this morphisms. It is the composition of elementary blocks which use square roots: 
\[\xymatrix{\relax
\mC \ar[d]_{f}       & \J \ar[d]_{f_{*}} \ar@{->>}[r] \ar@/^2pc/[rr]_{\psi} & \J/\left\{\pm1\right\} \ar[r] & \K \ar@{.>}[lld]  \\
\mE_{1}\times\mE_{2} & \mE_{1}\times\mE_{2}.
}\]

\bigskip

Let $P=\XYZT$ a point on the Kummer surface $\K_{\abcd}$ corresponding to the hyperelliptic curve $\mC$ of equation $y^2=f\left(x\right)$. We want the Mumford coordinates $\left(u,v\right)$ of the two opposite divisors $\psi^{-1}\left(P\right)=\left\{\pm D\right\}$. The divisors $D$ and $-D$ share the same polynomial $u$ and have opposite $v$ polynomial. First, compute the following quantities where for readability we note $\theta_{i}$ for $\theta_{i}\zero$:
\begin{eqnarray*}
\tse\z & = & \frac{-X\tci\thu+Y\tci\tdi-Z\tsi\tdi+T\tsi\thu}{\theta_{6}^{4}-\theta_{5}^{4}}, \\
\tne\z & = & \frac{ X\tsi\tdi-Y\tsi\thu+Z\tci\thu-T\tci\tdi}{\theta_{8}^{4}-\theta_{10}^{4}}, \\
\ton\z & = & \frac{ X\tsi\thu-Y\tsi\tdi+Z\tci\tdi-T\tci\thu}{\theta_{6}^{4}-\theta_{5}^{4}}, \\
\tdo\z & = & \frac{-X\tci\tdi+Y\tci\thu-Z\tsi\thu+T\tsi\tdi}{\theta_{8}^{4}-\theta_{10}^{4}}, \\
\ttz\z & = & \frac{-X\thu\tne+Y\tne\tdi+Z\tse\thu-T\tse\tdi}{\theta_{7}^{4}-\theta_{9}^{4}}, \\
\tqo\z & = & \frac{-X\tci\tne-Y\tsi\tse+Z\tci\tse+T\tsi\tne}{\theta_{7}^{4}-\theta_{9}^{4}}, \\
\tsz\z & = & \frac{-X\tsi\tse-Y\tci\tne+Z\tsi\tne+T\tci\tse}{\theta_{7}^{4}-\theta_{9}^{4}}. \\
\end{eqnarray*}
In the general case, $\tsz\z$ is non zero. In this case $u\left(\xi\right)=\xi^{2}+u_{1}\xi+u_{0}$ is of degree~$2$ and $v\left(\xi\right)=\pm\left(v_{1}\xi+v_{0}\right)$ is of degree $1$ with
\[ u_{0}=\la\frac{\thu \, \tqo\z}{\tdi \tsz\z}, \qquad u_{1}=\left(\la-1\right)\frac{\tci\, \ttz\z}{\tdi\tsz\z} -u_{0}-1 ,\]
\begin{eqnarray*}
v_{0}^{2} &=& -\frac{\theta_{1}^{4}\theta_{3}^{4}\thu\,\tqo\z }{\left(\tde\tqu\tdi\,\tsz\z\right)^{3}} \left(\tde\ttr\theta_{9}^{4} \, \tse\z \tdo\z+\tun\tqu\theta_{7}^{4}\,\tne\z\ton\z \right. \\
 & & \quad \left. +2\tun\tde\ttr\tqu\left(XZ+YT\right)  -\frac{\tun\ttr+\tde\tqu}{E'}\left(\left(X^2+Y^2+Z^2+T^2\right)- \right.\right.\\
 & & \quad \left. \left. -F\left(X T+Y Z\right)-G\left(XZ+YT\right)-H\left(XY+ZT\right)\right)\right).
\end{eqnarray*}
Since $u$ divides $v^2-f$ we obtain formul\ae{} to compute $v_{1}$. If $\tsz\z$ is zero and $\left(\la-1\right)\tci\,\ttz\z-\la\thu\,\tqo\z$ is non zero, then $u\left(\xi\right)=\xi+u_{0}$ is of degree $1$ and $v\left(\xi\right)=\pm v_{0}$ is constant with
\[ u_{0}=\frac{\la\thu\,\tqo\z}{\left(\la-1\right)\tci\,\ttz\z-\la\thu\,\tqo\z},\qquad v\left(\xi\right)=\pm v_{0}=\pm\sqrt{f\left(-u_{0}\right)}. \]
The last case is when $\tsz\z=\left(\la-1\right)\tci\,\ttz\z-\la\thu\,\tqo\z=0$. In that case, the divisor $D$ is the zero divisor.

\bigskip

Let $P$ be a point on the $\left(2,2\right)$-decomposable hyperelliptic curve $\mC$ given by the equation
\[ \mC: \quad \chi y^2=x\left(x-1\right)\left(x-\la\right)\left(x-\mu\right)\left(x-\nu\right), \quad \text{with } \la = \mu\frac{1-\nu}{1-\mu} . \]
We want to map $P$ to the elliptic curves. Let $\mC'$ be the curve given by
\[ \mC':\quad \kappa y^2=\left(x^2-1\right)\left(x^2-x_{2}^2\right)\left(x^2-x_{3}^2\right)  \]
where
\[ q=\pm\sqrt{\mu\left(\mu-\nu\right)},\  x_2= \frac{\mu+ q}{\mu- q}, \ x_{3}=\frac{1-\mu - q}{1-\mu + q}, \ \kappa= -\chi q\, \mu \left(\mu-1\right).\]
The curves $\mC$ and $\mC'$ are isomorphic by the change of coordinates
\[\begin{array}{lll}
 \mC              & \quad \longrightarrow \quad & \mC'    \\
 \left(x,y\right) & \quad \longmapsto     \quad & \left(\frac{x-\mu-q}{x-\mu+q},\frac{wy}{\left(x-\mu+q\right)^{3} }  \right)
\qquad \text{with} \ w=\frac{8q}{\left(\mu-q\right)\left(-1+\mu-q\right)}
\end{array}\]
The curve $\mC'$ maps to the elliptic curve $\mE: \ y^2=\left(x-1\right)\left(x-x_{2}^2\right)\left(x-x_{3}^2\right)$ by the morphism $\left(x,y\right)\mapsto\left(x^2,y\right)$. \\
Changing $q$ to $-q$ changes the elliptic curve $\mE$ to the other underlying elliptic curve. 

\bigskip

Let $f$ be the map from $\mC$ to the product of the two elliptic curves $\mE_{1}\times\mE_{2}$ given by the product of the two maps defined above. The push forward $f_{*}$ of $f$ is defined~by
\[ f_{*} \left\{
\begin{array}{lll}
 \J                                  & \quad \longrightarrow \quad & Jac\left(\mE_{1}\right)\times Jac\left(\mE_{2}\right)    \\
 D=\sum_{i=1}^{r} P_{i} -rP_{\infty} & \quad \longmapsto     \quad & \sum_{i=1}^{r} f\left(P_{i}\right) -r\,f\left(P_{\infty}\right)
\end{array}
\right.  \]
where we define $f\left(P_{\infty}\right)=\left(\mO_{\mE_{1}},\mO_{\mE_{2}}\right)$ to be the zero of the two elliptic curves. Note that the divisors in the Jacobian of the elliptic curves are not reduced. For elliptic curves, the Jacobian of the curve is isomorphic to the set of points on the curve. Thus we can identify $Jac\left(\mE_{i}\right)$ with $\mE_{i}$. We can rewrite $f_{*}$ as
\[ f_{*} \left\{
\begin{array}{lll}
 \J                                  & \quad \longrightarrow \quad & \mE_{1}\times \mE_{2}    \\
 D=\sum_{i=1}^{r} P_{i} -rP_{\infty} & \quad \longmapsto     \quad & \sum_{i=1}^{r} f\left(P_{i}\right) 
\end{array}
\right.  \]
A generic divisor $D$ in the Jacobian of a genus $2$ curve is the formal sum: $D=P_{1}+P_{2}-2P_{\infty}$. Thus in general we have to add the two points $f\left(P_{1}\right)$ and $f\left(P_{2}\right)$ on the elliptic curves. In practice, to do this operation, we need complete formul\ae{} for adding points on elliptic curves because we can't take a square root to get the $y$-coordinate which is needed for doubling in Weierstrass coordinates. Thus we use Jacobi or Edwards equations for the curves. Another solution to avoid square roots is to work in an extension ``field'' of degree $4$.

\bigskip

For stage $2$, GMP-ECM needs a point $(x_{2},y_{2})$ on the curve $y^2=x^3+Ax+B$ but we have a point $\left(x_{1}::z_{1}\right)$ on $\kappa zy^2=x^3+a_{2}x^2 z+a_{4}xz^2+a_{6}z^3$. First translate the point by $x\mapsto x-a_{2}z/3$ and divide the $x$-coordinate by $z$ to get a point $\left(x'_{1},\text{?}\right)$ on the curve $\kappa y^2=f\left(x\right)=x^3+a'_{4}x+a'_{6}$. Then the points $\left(x'_{1},\pm 1\right)$ are on the curve
\[ Dy^2=f\left(x\right)=x^3+a'_{4}x+a'_{6} \quad \text{with } D= f\left(x'_{1}\right). \]
By the change of variable 
\[ \left(x,y\right)\longmapsto \left(\frac{x}{D},\frac{y}{D}\right) \]
we get a point on the curve $y^2=x^3+Ax+B$ with $A=a'_{4}/D^2$ and $B=a_{6}/D^3$. Note that we can choose the sign of $y$ since a point and its opposite have the same order (remember that we cleared the power of $2$ in stage 1).

\bibliographystyle{amsplain}
\bibliography{article}

\ifx\NoAutoSpaceBeforeFDP\undefined\let\NoAutoSpaceBeforeFDP\relax\fi\def\noop%
sort#1{\errmessage{\string\noopsort\space is now forbidden !}}
  \ifx\dash\undefined\def\dash{\nobreak-\discretionary{}{}{}\hskip0pt}\fi
  \NoAutoSpaceBeforeFDP
  \ifx\bibfrench\undefined\def\biling#1#2{#1}\else\def\biling#1#2{#2}\fi
  \ifx\bibterse\undefined\def\rant#1#2{#1}\else\def\rant#1#2{#2}\fi
  \def\Inpreparation{\biling{In preparation}{en
  préparation}}\def\Preprint{\biling{Preprint}{préversion}}\def\Draft{\bilin%
g{Draft}{Manuscrit}}\def\Toappear{\biling{To appear}{À
  paraître}}\def\toappear{\biling{to appear}{à
  paraître}}\def\Toappearin{\biling{To appear in}{À paraître
  dans}}\def\Inpress{\biling{In press}{Sous presse}}\def\Seealso{\biling{See
  also}{Voir également}}\def\Maninprep{\biling{Manuscript in
  preparation}{Manuscrit en préparation}}\def\Availableat{\biling{Available
  at}{Disponible à l'adresse}}\def\availableat{\biling{available
  at}{disponible à l'adresse}}\def\Informalnote{\biling{Informal note}{Note
  informelle}}\def\researchreport{\biling{Research report}{Rapport de
  recherche}}\hyphenation{Ber-le-kamp}\providecommand\ieme{\up{ème}}
\providecommand{\bysame}{\leavevmode\hbox to3em{\hrulefill}\thinspace}
\providecommand{\MR}{\relax\ifhmode\unskip\space\fi MR }
\providecommand{\MRhref}[2]{%
  \href{http://www.ams.org/mathscinet-getitem?mr=#1}{#2}
}
\providecommand{\href}[2]{#2}
\begin{thebibliography}{10}

\bibitem{BeBiLaPe08}
D.~J. Bernstein, P.~Birkner, T.~Lange, and C.~Peters, \emph{{ECM} using
  {E}dwards curves}, Cryptology ePrint Archive, 2008,
  http://eprint.iacr.org/2008/016.

\bibitem{Duquesne07}
S.~Duquesne, \emph{Improving the arithmetic of elliptic curve in the {J}acobi
  model}, Inform. Process. Lett. \textbf{104} (2007), 101--105.

\bibitem{Gaudry07}
P.~Gaudry, \emph{Fast genus 2 arithmetic based on theta functions}, J. Math.
  Cryptol. \textbf{1} (2007), 243--265.

\bibitem{GaSc01}
P.~Gaudry and {\'E}.~Schost, \emph{On the invariants of the quotients of the
  {J}acobian of a curve of genus 2}, Applied Algebra, Algebraic Algorithms and
  Error-Correcting Codes (S.~Bozta\c{s} and I.~Shparlinski, eds.), Lecture
  Notes in Comput. Sci., vol. 2227, Springer-Verlag, 2001, pp.~373--386.

\bibitem{Kruppa08}
A.~Kruppa, \emph{Factoring into large primes with {P-1}, {P+1} and {ECM}},
  2008, Available at http://cado.gforge.inria.fr/workshop/slides/kruppa.pdf.

\bibitem{Lenstra87}
H.~W. {Lenstra, Jr.}, \emph{Factoring integers with elliptic curves}, Ann. of
  Math. (2) \textbf{126} (1987), 649--673.

\bibitem{Montgomery83}
P.~L. Montgomery, \emph{Evaluating recurrences of form
  $x_{m+n}=f\left(x_{m},x_{n},x_{m-n}\right)$ via {L}ucas chains}, 1983,
  Available at ftp.cwi.nl/pub/pmontgom/Lucas.ps.gz.

\bibitem{Montgomery87}
\bysame, \emph{Speeding the {Pollard} and {E}lliptic {C}urve {M}ethods of
  {F}actorization}, Math. Comp. \textbf{48} (1987), no.~177, 243--264.

\bibitem{Mumford83}
D.~Mumford, \emph{Tata lectures on theta {I}}, Progr. Math., vol.~28,
  Birkh\"{a}user, 1983.

\bibitem{Mumford84}
\bysame, \emph{Tata lectures on theta {II}}, Progr. Math., vol.~43,
  Birkh\"{a}user, 1984.

\bibitem{Shaska05}
T.~Shaska, \emph{Genus 2 curves covering elliptic curves, a computational
  approach}, Lecture Notes Ser. Comput. \textbf{13} (2005), 151--195.

\bibitem{Wamelen98}
P.~van Wamelen, \emph{Equations for the {J}acobian of a hyperelliptic curve},
  Trans. Amer. Math. Soc. \textbf{350} (1998), no.~8, 3083--3106.

\bibitem{ZiDo06}
P.~Zimmermann and B.~Dodson, \emph{20 years of {ECM}}, ANTS VII, Lecture Notes
  in Comput. Sci., vol. 4076, Springer-Verlag, 2006, pp.~525--542.

\end{thebibliography}

\end{document}